\documentclass[12pt]{amsart}
\usepackage{amssymb,amsfonts}
\usepackage{amsmath,amsthm}
\newtheorem{theorem}{Theorem}[section]
\newtheorem{proposition}[theorem]{Proposition}
\newtheorem{lemma}[theorem]{Lemma}
\newtheorem{corollary}[theorem]{Corollary}

\theoremstyle{definition}

\theoremstyle{remark}
\newtheorem{remark}[theorem]{Remark}

\numberwithin{equation}{section}

\def\La{{L^\infty_{\alpha+d}}}

\def\Ea{{E_{\alpha+d}}}

\def\R{\mathbb{R}}
\def\endProof{{\hfill$\Box$}}
\def\RR{{\mathcal{R}}}
\newcommand{\rf}[1]{{\rm (}\ref{#1}{\rm )}}

\begin{document}

\title[Convection equation with anomalous diffusion (\today)]{ Far field
 asymptotics of solutions\\
to convection equation with anomalous diffusion}

\author{Lorenzo Brandolese}
\address{Universit\'e de Lyon~; Universit\'e Lyon 1~;
CNRS UMR 5208 Institut Camille Jordan,
43 bd. du 11 novembre,
Villeurbanne Cedex F-69622, France.}
\email{brandolese{@}math.univ-lyon1.fr}
\urladdr{http://math.univ-lyon1.fr/$\sim$brandolese}

\author{Grzegorz Karch}
\address{Instytut Matematyczny, Uniwersytet Wroc\l awski pl.
Grunwaldzki 2/4, 50-384 Wroc\l aw, Poland}
\email{karch{@}math.uni.wroc.pl}
\urladdr{http://www.math.uni.wroc.pl/$\sim$karch}

\thanks{
The preparation of this paper  
was supported in part by the European Commission Marie Curie Host
Fellowship
for the Transfer of Knowledge ``Harmonic Analysis, Nonlinear
Analysis and Probability''  MTKD-CT-2004-013389,
and in part by the Polonium Project PAI EGIDE N.~09361TG.
The first author gratefully thanks the Mathematical Institut of Wroc\l
aw University
for the warm hospitality. The preparation of this paper by the second
author
was also partially
supported by the
grant N201 022 32 / 09 02.
}

\subjclass[2000]{Primary 35K; Secondary 35B40, 35Q, 60H}

\date{\today}

\keywords{Anomalous diffusion,
asymptotic profiles, self-similar solutions, decay estimates, fractal
Burgers equation,
conservation laws.}

\begin{abstract}
The initial value problem for the conservation law
$\partial_t u+(-\Delta)^{\alpha/2}u+\nabla \cdot f(u)=0$
is studied for $\alpha\in (1,2)$ and under natural polynomial growth
conditions imposed on the nonlinearity. We find
the asymptotic expansion as $|x|\to \infty$ of solutions to this
equation corresponding to initial conditions, decaying sufficiently
fast at infinity.
\end{abstract}

\maketitle
%%%%%%%%%%%%%%%%%%%%%%%%%%%%%%%
%%%%%%%%%%%%%%%%%%%%%%%%%%%%%%%
%%%%%%%%%%%%%%%%%%%%%%%%%%%%%%%

\section{Introduction}
We study properties of solutions to the initial value problem for the 
multidimensional conservation law  with the anomalous diffusion
\begin{eqnarray}
\label{eq}
&&\partial_t u+(-\Delta)^{\alpha/2}u+\nabla \cdot f(u)=0, \quad x\in
\R^d, t>0,\\
\label{ini}
&&u(x,0)=u_0.
\end{eqnarray}
Here, we always impose the standing assumption
$1<\alpha<2$. Moreover, we assume that the $C^1$- vector field
$f(u)=\bigl(f_1(u),\ldots,
f_d(u)\bigr)$ is of a polynomial
growth, namely, it satisfies the usual estimates
\begin{eqnarray}
\label{hf1}
|f(u)|\le C|u|^q
\quad\mbox{and}\quad
\label{nla}
|f(u)-f(v)|\le C|u-v|\bigl(|u|^{q-1}+|v|^{q-1}\bigr)
\end{eqnarray}
for some constants $C>0$, $q>1$ and for all $u,v\in\R$
(in fact,  assumption \rf{hf1} can be slightly relaxed   in some parts
of our considerations, cf. Remark \ref{2.3}, below).

Linear evolution problems involving fractional Laplacian
describing {\it  the ano\-ma\-lous diffusion} (or $\alpha$-stable L\'evy
diffusion)
have been extensively
studied in the mathematical and physical literature (see, e.g.,
\cite{MK00}).
The probabilistic interpretation of nonlinear
evolution problems with an anomalous diffusion, obtained recently by
Jourdain, M\'{e}l\'{e}ard, and Woyczy\'{n}ski \cite{JMW05}, motivated us
to study \rf{eq}-\rf{ini}.
The authors of \cite{JMW05} considered a class of nonlinear
integro-differential equations involving a
fractional power of the Laplacian and a nonlocal quadratic
nonlinearity  represented by a singular integral
operator. They associated with the equation a nonlinear
singular diffusion
and proved {\it propagation
of chaos} to the law of this diffusion for the related 
interacting particle systems.
In particular, due to the probabilistic
origin of \rf{eq}-\rf{ini}, 
the
function $u(\cdot,t)$ should be
interpreted as  the density of a probability
distribution for every
$t>0$,
if the initial datum is so.

Although, the motivation for this paper comes from the probability
theory, our methods are purely analytic.
Hence, if $X(t)$ is the symmetric $\alpha$-stable L\'evy process, its
density of the probability distribution  $p_\alpha(x,t)$
is the fundamental solution
of the linear equation
\begin{equation}
\label{leq}
 \partial_t v+(-\Delta)^{\alpha/2}v=0,
\end{equation}
hence, $p_\alpha$
can be computed {\it via} the Fourier
transform
$ \widehat p_\alpha(\xi,t)=e^{-t|\xi|^\alpha}$.
In particular,
\begin{equation*}
p_\alpha(x,t)=t^{-d/\alpha}P_\alpha(xt^{-1/\alpha}),
\end{equation*}
where $P_\alpha$ is the inverse Fourier transform of $e^{-
|\xi|^\alpha}$ (see \cite[Ch.~3]{J1} for more details).
It is well known that for every $\alpha\in (0,2)$ the function
$P_\alpha$
is  smooth, nonnegative, and
satisfies the estimates
\begin{equation}
\label{EP}
0<P_\alpha(x)\le C(1+|x|)^{-(\alpha+d)} \qquad \hbox{and}\qquad
|\nabla P_\alpha(x)|\le C(1+|x|)^{-(\alpha+d+1)}
\end{equation}
for a constant $C$ and all $x\in\R^d$.
Moreover,
\begin{equation}\label{as P}
P_\alpha(x)=c_0|x|^{-(\alpha+d)}+O\left(|x|^{-(2\alpha+d)}\right),
\qquad\hbox{as $|x|\to\infty$},
\end{equation}
and
\begin{equation}\label{as nab P}
\nabla
P_\alpha(x)=-c_1\,x|x|^{-(\alpha+d+2)}+O\left(|x|^{-(2\alpha+d+1)}
\right), \qquad\hbox{as $|x|\to\infty$,}
\end{equation}
where
$$ c_0=\alpha2^{\alpha-1}\pi^{-(d+2)/2}\sin
(\alpha\pi/2)\Gamma\Bigl(\frac{\alpha+d}{2}\Bigr)
\Gamma\Bigl(\frac{\alpha}{2}\Bigr),$$
and
$$ c_1=2\pi\alpha2^{\alpha-1}\pi^{-(d+4)/2}\sin
(\alpha\pi/2)\Gamma\Bigl(\frac{\alpha+d+2}{2}\Bigr)
   \Gamma\Bigl(\frac{\alpha}{2}\Bigr).$$
We refer to~\cite{BluG60} for a proof of the  formula
(\ref{as P}) with the explicit constant $c_0$.
The optimality of the estimate of the lower
order term in~\rf{as P} is due Kolokoltsov \cite[Eq. (2.13)]{Kol99},
where  higher order expansions of $P_\alpha$
are also computed.
The proof of the asymptotic expression (\ref{as nab P}) and
the value of $c_1$
can be  deduced from (\ref{as P}) using an identity by 
Bogdan and Jakubowski
\cite[Eq. (11)]{Jak06}.

The asymptotic formula (\ref{as P}) for the kernel  $P_\alpha$ plays an
important role  in the theory of $\alpha$-stable
processes.
The main goal in this
work is to present a method which allows to derive 
analogous asymptotic
expansions
as
$|x|\to
\infty$
of solutions  to the Cauchy problem \rf{eq}-\rf{ini}.
In the next section, we recall several properties of solutions to
\rf{eq}-\rf{ini} and we state our main results: Theorems \ref{th21}
and
\ref{th22}. In Section \ref{space-time},
we gather  technical
space-time
estimates of solutions to \rf{eq}-\rf{ini}. The proofs of
Theorems \ref{th21} and \ref{th22} are contained in Section
\ref{secP}.

{\bf Notation.}
The $L^p$-norm   of a Lebesgue
measurable, real-valued function $v$ defined on $\R^d$
is denoted by $\|v\|_p$.
In the following, we use the weighted $L^\infty$ space
\begin{equation}\label{La}
L^\infty_{\vartheta}= \{v\in L^\infty(\R^d)\;:\;
\|v\|_{L^\infty_\vartheta}\equiv
{\rm ess\, sup}_{x\in\R^d}|v(x)|(1+|x|)^{\vartheta}<\infty\},
\end{equation}
for any $\vartheta\ge0$,
and its homogeneous counterpart
\begin{equation*}
\dot L^\infty_{\vartheta}= \{v\in L^\infty_{\rm loc}(\R^d\backslash\{0\})\;:\;
\|v\|_{\dot L^\infty_\vartheta}\equiv
{\rm ess\, sup}_{x\in\R^d}|v(x)||x|^{\vartheta}<\infty\}.
\end{equation*}
The constants (always independent
of $x$) will be
denoted by the
same letter $C$, even if they may vary from line to line.
Sometimes, we write,  e.g.,  $C=C(T)$ when we want to
emphasize
the dependence of $C$ on a parameter $T$.

%%%%%%%%%%%%%%%%%%%%%%%%%%%%%%%
\section{Main results}

It is well known (see \cite{BKWstudia,DGV,DroI06}) that given
$u_0\in L^1(\R^d)$ and  $1<\alpha\le2$,  the initial value
problem~\eqref{eq}-\eqref{ini} has the unique
solution $u\in C([0,\infty), L^1(\R^d))$.
Moreover, this solution satisfies
\begin{equation*}
u\in C((0,\infty), W^{1,p}(\R^d))
\end{equation*}
for every $p\in [1,\infty]$ and the following estimates hold true
(see \cite[Cor. 3.2]{BKWstudia})
\begin{equation}
\label{pbo}
\| u(t)\|_p \le C t^{-\frac{d}{\alpha}(1-\frac{1}{p})}\|u_0\|_1
\end{equation}
for all $t>0$ and $C$ independent of $t$ and of $u_0$.
Under the additional assumption $u_0\in L^p(\R^d)$, the corresponding
solution satisfies $u\in C([0,\infty), L^p(\R^d))$ together with the
estimate
\begin{equation}
\label{Lp bound}
\| u(t)\|_p\le \|u_0\|_p.
\end{equation}

Below, in Proposition \ref{prop1},
we complete these preliminary results
providing the estimates of solutions to \rf{eq}-\rf{ini} in weighted
$L^\infty$-spaces.
In particular, if $u_0\in\La$ (cf. \rf{La}),
then
the corresponding solution of \rf{eq}-\rf{ini} satisfies
$u\in C([0,T], \La)$
for every $T>0$.
Such a result was already obtained in the one dimensional case, see
\cite[Sect. 2]{HKNS05}.
In Section \ref{space-time}, we state and prove its multidimensional
counterpart for the
completeness of the exposition.
We complement this result with additional estimates for the
gradient
of the solution, which will be useful in the proofs of asymptotic
formulas in  
Section~\ref{secP}.

Let us recall that, when studying the large time behavior
of solutions for the problem~\rf{eq}-\rf{ini},
an important role is played by the {\it critical exponent}
$$\widetilde{q}\equiv 1+\frac{\alpha-1}{d}.$$
Indeed, using the terminology of \cite{BKW01}
the behavior of solutions as $t\to\infty$ is
\emph{genuinely non-linear\/} when
$q=\widetilde{q}$, 
is \emph{weakly non-linear\/} when $q>\widetilde{q}$,
and is (expected to be) \emph{hyperbolic\/} when $1<q<\widetilde{q}$.

In this paper, in the supercritical case $q>\widetilde{q}$, as well as
for $q=\widetilde{q}$
provided $\|u_0\|_1$ is sufficiently small, we
will improve the space-time estimates
of \cite[Sect. 2]{HKNS05},
showing that
\begin{equation}
|u(x,t)|\leq Cp_\alpha(x,1+t)\label{u p a},
\end{equation}
for all $x\in\R^d$, $t>0$, and $C>0$ independent of $x,t$.
Under the additional assumption
that $\nabla u_0\in L^\infty_{\alpha+d+1}$, we will also prove that
\begin{equation}
\label{nupi}
\|\nabla u(t)\|_{L^\infty_{\alpha+d+1}}
\le C (1+t),
\end{equation}
see Theorem \ref{th up}, below.
In other words, $\nabla u(x,t)$ has the same space-time
decay profiles as $\nabla p_\alpha(x,1+t)$ (cf. the second inequality
in \rf{EP}).

Furthermore,
we make evidence of the {\it second critical exponent},
namely,
$$q^*\equiv 1+ \frac{1}{\alpha+d},$$
playing an important role in the study of the pointwise behavior of
solutions as $|x|\to\infty$.
The following theorem explains the role of $q^*$,
showing that any decaying solution has a precise spatial
asymptotic profile.
Here, we denote by
$S_\alpha(t)u_0(x)=p_\alpha(t)*u_0(x)$ the solution of the linear
equation~\rf{leq} supplemented with the initial datum~$u_0$.

\begin{theorem}
\label{th21}
Assume that $\alpha\in (1,2)$,
and let $u= u(x,t)$ be the  solution of~\eqref{eq}-\eqref{ini}
with the nonlinearity satisfying \rf{hf1},
and with $u_0\in \La$.

(i)
Then, for all $t>0$, $x\in\R^d$,
\begin{equation}
\label{pr u 2}
\begin{split}
u(x,t)= S_\alpha(t)u_0(x)+ &\frac{c_1x}{|x|^{\alpha+d+2}}\cdot
 \int_0^t \!\! \int (t-s) f(u(y,s))\,dy\,ds \\
+&O\Bigl( \max\bigl\{ |x|^{-q(\alpha+d)} \,;\,|x|^{-(\alpha+d+2)}
\bigr\} \Bigr),
\qquad\hbox{as $|x|\to\infty$},
\end{split}
\end{equation}
 uniformly in
any time interval $t\in [0,T]$, $T>0$.
This conclusion is interesting only when the last term on the right
hand side
of~\rf{pr u 2}
is the lower order term as $|x|\to\infty$: this happens when
$q>q^*$.

(ii)
The conclusion at the point (i) can be improved
under the additional assumption
$\nabla u_0\in L_{\alpha+d+1}^\infty$,
replacing the third term on the right hand side of~\rf{pr u 2}
by
\begin{equation*}
O\Bigl( \max\bigl\{ |x|^{-q(\alpha+d)-1} \,;\,|x|^{-(\alpha+d+2)}
\bigr\} \Bigr),
\qquad\hbox{as $|x|\to\infty$}.
\end{equation*}
Now, this conclusion is interesting also for
 $1< q \le q^*$.

(iii) If
$u$ satisfies inequality \rf{u p a} for all $x\in \R^d$, $t>0$, and
$C>0$ independent of
$x,t$,
then one can make precise the behavior for large~$t$ of the remainder term
in
relation~\rf{pr u 2},
replacing it by
$$O\Bigl((1+t)^N\max\bigl\{|x|^{-q(\alpha+d)}\,;\, |x|^{-(\alpha+d+2)}
\bigr\}\Bigr),
\quad\hbox{as $|x|\to\infty$},$$
uniformly in $t\in[0,\infty)$,
for some exponent $N=N(\alpha,q,d)\le 3$, independent on $u_0$.

If, in addition, the solution satisfies inequality~\rf{nupi},
 the conclusion at the point (ii) can be improved
replacing the remainder term by
$$O\Bigl((1+t)^N\max\bigl\{|x|^{-q(\alpha+d)-1}\,;\, |x|^{-
(\alpha+d+2)} \bigr\}\Bigr)$$
for some exponent $N=N(\alpha,q,d)\le 3$, independent on $u_0$, and
the convergence as $|x|\to\infty$ holds true uniformly in $t\in
[0,\infty)$.
\end{theorem}

It follows from the Duhamel formula that the solution of the Cauchy
problem \rf{eq}-\rf{ini}
satisfies the  integral equation
\begin{equation}
\label{duham}
u(t)=S_\alpha(t)u_0-\int_0^t \nabla S_\alpha(t-s)\cdot f(u)(s)\,ds.
\end{equation}
Hence, it is possible to give a heuristic explanation of the role of
space-critical
exponent $q=q^*$, simply, by  looking at the integrand
of the second term on the right hand side of~\rf{duham}.
Indeed, the kernel of~$\nabla S_\alpha$ behaves as $|x|^{-
(\alpha+d+1)}$
as $|x|\to\infty$ (cf. \rf{as P}), whereas
$|f(u(x,t))|\leq C |x|^{- q(d+\alpha)}$ for $u(t)\in \La$.
Then, it is natural to expect that the large space asymptotics
is influenced by the competition between  these two
decay rates as $|x|\to\infty$.
In fact,
 the proof of Theorem \ref{th21} (given in Section \ref{secP})
consists in finding  the asymptotic expansion of the
second term on the right hand side of \rf{duham} and 
the equality between these two decay rates
occurs precisely when $q=q^*$.

\begin{remark}
It is worth observing that this type of asymptotic expansion of
solutions to convection-diffusion equations
 is specific of  the fractional nature
of the diffusion operator $(-\Delta)^{\alpha/2}$
and is caused by the algebraic decay of the fundamental solution
$p_\alpha(x,t)$.
For the viscous Burgers equation, or for multidimensional
diffusion-convection equations with standard dissipation
({\it i.e.\/}, with the usual Laplacian)
Theorem~\ref{th21} remains valid, but it is not interesting because
the coefficient~$c_1$ vanishes in the limit case $\alpha=2$.
  \end{remark}

\begin{remark}\label{2.3}
The conclusion~(i) of Theorem~\ref{th21} remains valid under more general assumptions
on the nonlinearity. What we really need is that $f$ is a $C^1$-vector field
such that $|f(u)|\le c(R)|u|^q$ for some  $q>1$, a
continuous nondecreasing function $c(\cdot)$ on $[0,\infty)$, and all $|u|\le R$.
For the part~(ii), we need also a similar condition for $f'$, namely,  $|f'(u)|\le c_1(R)|u|^{q-1}$
for $|u|\le R$. 
On the other hand, the present form of  Theorem \ref{th21}.iii is 
 no longer valid for such more general nonlinearities.
%The growth of the remainder terms as $t\to\infty$ being dependent on~$f$, conclusions~(iii)-(iv)
%are no longer valid, in that form, in this more general situation.
Our more stringent assumption~\eqref{hf1} allows us
to present the essential ideas avoiding uninteresting technicalities in the proofs, in
particular, separating the cases of large and small~$u$ in our estimates.
Moreover, such an assumption is well suited for studying self-similar solutions.
\end{remark}

For the homogeneous nonlinear term $\nabla \cdot f(u)= b\cdot
\nabla(u|u|^{q-1})$ with
a fixed $b\in \R^d$ and with the time-critical exponent
$q=\widetilde{q}$, the authors
 of \cite{BKW01}
constructed a family of self-similar solutions $u_M=u_M(x,t)$ of
equation \rf{eq}.
Those functions satisfy the scaling relation
\begin{equation}\label{uM self}
u_M(x,t)=t^{d/\alpha} U_M(xt^{-1/\alpha})\quad \mbox{where} \quad
U_M(x)=u_M(x,1)
\end{equation}
for all  $x\in\R^d$ and  $t>0$. Moreover,
 each of them  is the unique solution of the initial value problem
\begin{eqnarray}
&&\partial_t u +(-\Delta)^{\alpha/2}u+b\cdot
\nabla(u|u|^{(\alpha-1)/d})=0\label{eq-sing}\\
&&u(x,0)=M\delta_0\label{ini-sing}
\end{eqnarray}
for $\alpha\in (1,2)$ and $M>0$, where $\delta_0$ is the Dirac delta.
We refer the reader to \cite{BKW01} for more information concerning
solutions of
problem \rf{eq-sing}-\rf{ini-sing}.

In this paper, 
we complete  results from \cite{BKW01} providing  space-time estimates
of those self-similar solutions.
First, in Corollary \ref{cor est U} below, we establish,
for sufficiently small $M>0$, the estimate
\begin{equation}\label{self u p 0}
0\leq u_M(x,t)\leq Cp_\alpha(x,t) \quad\hbox{for all $x\in\R^d$ and
$t>0$},
\end{equation}
We conjecture that such estimate remains true without the smallness
assumption
imposed on $M$.
Inequality \rf{self u p 0}
plays a crucial role in the proof of the following  asymptotic
expansion
of the self-similar kernel $U_M$.

\begin{theorem} \label{th22}
Assume that $1<\alpha<2$ and $\widetilde{q}>q^*$.
Let $u_M$ be a self-similar solution of~\rf{eq-sing}-\rf{ini-sing},
satisfying the estimate~\rf{self u p 0}.
Then the self-similar profile $U_M(x)=u_M(x,1)$
has the following behavior as $|x|\to\infty$:
\begin{equation}
\label{assud}
U_M(x)=
MP_\alpha(x) +\frac{c_1\alpha^2}{\alpha+1} \,
\|U_M\|_{\widetilde{q}}^{\widetilde{q}} \,
\frac{b\cdot x}{|x|^{\alpha+d+2}}
+O\Bigl( \max\bigl\{ |x|^{-\widetilde{q}(\alpha+d)} \,;\,
	|x|^{-(\alpha+d+2)} \bigr\}
	\Bigr).
\end{equation}
\end{theorem}

The asymptotic expansion of solutions to \rf{eq} stated in \rf{pr u 2} 
and in \rf{assud} can be viewed as the true counterparts of the 
well-known result for the $\alpha$-stable distribution recalled in 
\rf{as P}.

%%%%%%%%%%%%%%%%%%%%%%%%%%%%%%%

\section{Preliminary space-time estimates}
\label{space-time}

We begin this section by the study of
  the solution of the linear problem
\begin{equation}
\partial_t v+(-\Delta)^{\alpha/2}v=0, \qquad v(x,0)=v_0
\end{equation}
denoted by
$$ v(x,t)=S_\alpha(t)v_0(x)=p_\alpha(\cdot,t)*v_0(x).$$
The following lemma contains a direct
generalization to $\R^d$ of estimates from
\cite[Lemma 1.40]{HKNS05}.
By this reason, we sketch its proof only.

\begin{lemma}
\label{LL1}
Assume that $v_0\in L^\infty_{\alpha+d}$.
There exists $C>0$ independent of $v_0$ and $t$ such that
\begin{eqnarray}
\label{St 1}
\|S_\alpha(t)v_0\|_\infty &\le& C\min \left\{ t^{-d/\alpha}\|v_0\|_1,
\|v_0\|_\infty\right\},\\
\label{St 2}
\|S_\alpha(t)v_0\|_{L^\infty_{\alpha+d}}&\le&
C(1+t)\|v_0\|_{L^\infty_{\alpha+d}} \,,\\
\label{St 3}
\|\nabla S_\alpha(t) v_0\|_{L^\infty_{\alpha+d}} &\le&
Ct^{-1/\alpha}\|v_0\|_{L^\infty_{\alpha+d}}
    +Ct^{1-1/\alpha}\|v_0\|_1 \,,
\end{eqnarray}
\end{lemma}

\proof
Estimate \rf{St 1} results immediately from the Young
inequality
applied  to the convolution $S_\alpha(t)v_0=p_\alpha(t)*v_0$, due
to the identities
$$
\|p_\alpha(t)\|_1=1, \quad  \|p_\alpha
(t)\|_\infty = t^{-d/\alpha}  \|P_\alpha\|_\infty
\quad \mbox{for all $t>0$.}
$$

Since
$|v_0(x)|\leq C(1+|x|)^{-
(d+\alpha)}$, by the asymptotic properties of the kernel
$p_\alpha(x,1)=P_\alpha(x)$ (cf. \rf{as P}), we immediately obtain
$|v_0(x)|\leq
Cp_\alpha(x,1)$ for  all $x\in \R^d$ and a constant $C>0$ independent
of $x$.
Consequently, by the semigroup property, we conclude
$$
\|S_\alpha(t)v_0\|_\La\leq C \|S_\alpha (t)p_\alpha(1)\|_\La=
C\|p_\alpha(t+1)\|_\La \leq C(1+t).
$$
Now, replacing $v_0$ by
$v_0/\|v_0\|_\La$ we obtain \rf{St 2}.

To prove \rf{St 3}, we use the pointwise estimate
$$
(1+|x|)^{\alpha+d}\leq C(1+|y|)^{\alpha+d}+C |x-y|^{\alpha+d},
$$
valid for  all $x,y\in \R$ and a constant $C>0$,
and we apply the Young inequality.
We get
$$
\|\nabla S_\alpha(t) v_0\|_{L^\infty_{\alpha+d}} \le
C \|\nabla p_\alpha(t)\|_1  \|v_0\|_{L^\infty_{\alpha+d}}
    +C \|\nabla p_\alpha(t)\|_{\dot L^\infty_{\alpha+d}} \|v_0\|_1
    $$
and \rf{St 3} immediately follows.   
\endProof

\medskip
Under an additional information on the gradient of $v_0$,
we can obtain analogous estimates for $\nabla S_\alpha(t)v_0$.
In order to give a precise statement, let us introduce the space
\begin{equation}\label{Ea}
\Ea\equiv \{v\in W^{1,\infty}_{loc}(\R^d)\;:\; \|v\|_{\Ea}\equiv
\|v\|_{\La}
+\|\nabla v\|_{L^\infty_{\alpha+d+1}}
<\infty\}.
\end{equation}

\begin{lemma}
\label{LLL1}
Assume that $v_0\in \Ea$.
There exists $C>0$ independent of $v_0$ and $t$ such that
\begin{eqnarray}
\label{ST1}
\|\nabla S_\alpha(t)v_0\|_\infty &\le& C\min \left\{
t^{-(d+1)/\alpha}\|v_0\|_1\,;\,
t^{-1/\alpha}\|v_0\|_\infty\,;\,
\|\nabla v_0\|_\infty\right\},\\
\label{ST2}
\|S_\alpha(t)v_0\|_\Ea &\le& C(1+t)\|v_0\|_\Ea \,,\\
\label{ST3}
\|\nabla S_\alpha(t) v_0\|_\Ea &\le&
 Ct^{-1/\alpha}\|v_0\|_\Ea    +Ct^{1-1/\alpha}\|v_0\|_1\,
\end{eqnarray}
for all $t>0$.
\end{lemma}

\proof
Estimate~\rf{ST1} is the straightforward application of
the $L^1$-$L^\infty$ convolution inequalities.
In order to prove \rf{ST2}
using the radial symmetry of $p_\alpha(\cdot,t)$,  we see that, for
all $R>0$,
 $\int_{B_R} \nabla p(y,t)\,dy=0$,  where $B_R$ denotes the ball
centered at the origin and  of radius $R$.
Hence,
$$\nabla S_\alpha(t)v_0(x)=\int_{|y|\le |x|/2} \bigl[
v_0(x-y)-v_0(x)\bigr]\nabla p_\alpha(y,t)\,dy
 +\int_{|y|\ge |x|/2}  v_0(x-y)\nabla p_\alpha(y,t)\,dy.$$
This decomposition shows that, for some constant $C>0$, the quantity
$|\nabla S_\alpha(t)v_0(x)|$ can be bounded from above by
\begin{equation*}
\begin{split}
C|x|^{-(\alpha+d+1)} \|\nabla v_0\|_{L^\infty_{\alpha+d+1}}
\int_{\R^d}|y|\,|\nabla p(y,t)|\,dy
+Ct|x|^{-(\alpha+d+1)}\int_{\R^d}|v_0(y)|\,dy,
\end{split}
\end{equation*}
which implies
\begin{equation}
\label{st2c}
\| \nabla S_\alpha(t)v_0\|_{\dot L^\infty_{\alpha+d+1}}
 \le
 C\left( \|\nabla v_0\|_{\dot L^\infty_{\alpha+d+1}}+
t\|v_0\|_1\right)
 \le
 C(1+t)\|v_0\|_{\Ea}.
\end{equation}
Now, estimate~\rf{ST2} follows from~\rf{St 2}, \rf{ST1} and from the
bound
for the homogeneous norm \rf{st2c}.

Let us prove~\rf{ST3}. By~\rf{St 3}
and the inequality
$$ \|\nabla^2S_\alpha(t)v_0\|_{\infty}
\le \| \nabla S_\alpha(t)\|_1\| \nabla v_0\|_\infty
\le Ct^{-1/\alpha} \| \nabla v_0\|_\infty,$$
 we see that
we only have to establish the following estimate in the homogeneous 
space
$\dot L^\infty_{\alpha+d+1}$
\begin{equation}
\label{st3c}
\|\nabla^2S_\alpha(t)v_0\|_{\dot L^\infty_{\alpha+d+1}}
\le Ct^{-1/\alpha} \|v_0\|_\Ea +Ct^{1-1/\alpha}\|v_0\|_1.
\end{equation}
To prove~\rf{st3c}, we consider the decomposition
$$\nabla^2 S_\alpha(t)v_0(x)=(J_1+J_2+J_3)(x,t),$$
where
\begin{equation*}
\begin{split}
J_1(x,t)&\equiv\int_{|y|\le |x|/2}
[v_0(x-y)-v_0(x)]\nabla^2p_\alpha(y,t)\,dy,\\
J_2(x,t)&\equiv\int_{|y|\ge |x|/2}
v_0(x-y)\nabla^2p_\alpha(y,t)\,dy,\\
J_3(x,t)&\equiv -v_0(x)\int_{|y|\ge |x|/2}\nabla^2 p_\alpha(y,t)\,dy
\end{split}
\end{equation*}
(note that $\int_{\R^d}\nabla^2 p_\alpha(y,t)\;dy=0$).
>From the well known estimate (see \cite{Kol99})
\begin{equation}
\label{n2p}
|\nabla^2 P_\alpha(x)|\le C(1+|x|)^{-(\alpha+d+2)},
\end{equation}
we deduce
$\int_{\R^d} |y|\, |\nabla^2 p_\alpha(y,t)|\,dy\le C t^{-1/\alpha}$.
Then, the application of the Taylor formula in the integral defining
$J_1$
yields
$$ |J_1(x,t)|\le Ct^{-1/\alpha}|x|^{-(\alpha+d+1)}\|\nabla
v_0\|_{L^\infty_{\alpha+d+1}}.$$

To deal with the terms  $J_2$ and $J_3$, we use two different
pointwise estimates
of $\nabla^2 p_\alpha(x,t)$ resulting from~\rf{n2p}:
\begin{equation*}
|\nabla^2 p_\alpha(x,t)|\leq 
Ct^{-(d+2)/\alpha}\left(1+|x|t^{-1/\alpha}\right)^{-(\alpha+d+2)}\leq
Ct^{1- 1/\alpha}|x|^{-(\alpha+d+1)}
\end{equation*}
and
\begin{equation*}
|\nabla^2 p_\alpha(x,t)|\leq  
Ct^{- 1/\alpha}|x|^{-(d+1)},
\end{equation*}
which imply
\begin{equation*}
\begin{split}
|J_2(x,t)|
&\leq \sup_{|y|\geq |x|/2} |\nabla^2p_\alpha(y,t)| 
\int_{|y|\geq |x|/2}|v_0(x-y)|\;dy\\
&\le Ct^{1-1/\alpha} |x|^{-(\alpha+d+1)} \|v_0\|_1
\end{split}
\end{equation*}
and
\begin{equation*}
\begin{split}
|J_3(x,t)|&\leq C|x|^{-(\alpha+d)}\|v_0\|_\La
\int_{|y|\ge |x|/2}|\nabla^2 p_\alpha(y,t)|\,dy\\
&\le Ct^{-1/\alpha}|x|^{-(\alpha+d+1)}\|v_0\|_{L^\infty_{\alpha+d}}.
\end{split}
\end{equation*}
Combining all these inequalities yields~\rf{st3c}.
\endProof

\medskip

We are in a position to construct solutions of the Cauchy problem
\rf{eq}-\rf{ini} in the weighted space $\La$.

\begin{proposition}
\label{prop1}
(i)
Let $\alpha\in(1,2)$ and $q>1$.
Assume that $u$ is a solution of the Cauchy problem~\eqref{eq}-
\eqref{ini} with the nonlinearity
satisfying~\eqref{hf1}.
If $u_0\in L^\infty_{\alpha+d}$, then
\begin{equation}
\label{d sol}
u\in C([0,T], L^\infty_{\alpha+d}) \quad\hbox{for each $T>0$}.
\end{equation}

(ii)
Under the more stringent assumption $u_0\in E_{\alpha+d}$, cf. 
\rf{Ea},
 we have also
 \begin{equation}
\label{d grad}
 u\in L^\infty([0,T], \Ea)
 \quad\hbox{for each $T>0$}.
\end{equation}
\end{proposition}

\proof
 In order to prove
\rf{d sol},
it suffices
to show that the nonlinear operator
$$
T(u)(t)= S_\alpha (t)u_0-\int_0^t \nabla S_\alpha(t-\tau)
f(u(\tau))\;d\tau
$$
has the fixed point in the space
$$
X_T=\{u\in C([0,T],\La)\;:\; \sup_{t\in[0,T]}\|u(t)\|_\La<
\infty\}.
$$
As usual, we work in the ball
$B(0,R)=\{u\in C([0,T],\La)\;:\; \sup_{t\in[0,T]}\|u(t)\|_\La\leq
R\}$, 
where $R=M\|u_0\|_{L^\infty_{\alpha+d}}$ and $M>0$ is a large constant, and $T>0$.
Combining inequality~\rf{St 3} with assumption~\rf{nla}
we get
\begin{eqnarray}
 \label{ase}
\nonumber
\|\nabla S_\alpha(t) f(u)\|_{\La} &\leq&
 Ct^{-1/\alpha}\|\,|u|^q\|_{\La}+Ct^{1-1/\alpha}\|u\|_q^q\\
&\leq&
Ct^{-1/\alpha}(1+t)\|u\|_\infty^{q-1}\|u\|_{\La}.
\end{eqnarray}
Applying now inequalities \rf{St 2}-\rf{ase}
we can estimate, for $u\in B(0,R)$,
\begin{eqnarray*}
\|T(u)(t)\|_{\La}&\leq& C(1+t)\|u_0\|_{\La} \\
&&\quad+ C R^{q-1} \int_0^t (t-\tau)^{-1/\alpha}(1+(t-\tau)) \|u(\tau)\|_{\La}\;d\tau\\
&\leq&
R/2+CM^{q-1}\|u_0\|_{L^\infty_{\alpha+d}}^{q-1} R\,t^{1-1/\alpha}(1+t)\\
&\leq& R,
\end{eqnarray*}
provided that $0\le t\le T$ and
$$T\le C\min\{1,\|u_0\|_{L^\infty_{\alpha+d}}^{-\alpha(q-1)/(\alpha-1)} \},$$
with $C>0$ small enough.

In the same way,  for all $u,\tilde u\in B(0,R)$,
\begin{eqnarray*}
\|T(u)(t)-T(\tilde u)(t)\|_{\La}\leq
 CR^{q-1} \int_0^t (t-\tau)^{-1/\alpha} (1+(t-\tau))
\|u(\tau)-\tilde u(\tau)\|_{\La}\;d\tau.
\end{eqnarray*}
The Banach fixed point theorem now guarantees the existence of a local-in-time solution.
In the next step,
such solution must be extended globally-in-time.
The argument is standard: we fix $T>0$ arbitrarily large
and using that $\|u(t)\|_\infty \le C$ on $[0,T]$ (see inequality~\eqref{Lp bound}),
we show that $\|u(t)\|_{L^\infty_{\alpha+d}}$ does not blow up on $[0,T]$.
Indeed for some constants $C_1$, $C_2$, \dots, depending on~$T$, for $0\le t\le T$ we have
$$ \|u(t)\|_{\La}\le C_1+C_2\int_0^t (t-\tau)^{-1/\alpha}\|u(\tau)\|_{\La}\,d\tau.$$
Iterating this inequality and applying Fubini's theorem we get
$$ \|u(t)\|_{\La}\le C_3+C_4\int_0^t (t-\tau)^{1-2/\alpha}\|u(\tau)\|_{\La}\,d\tau.$$
We  repeat this argument until 
we obtain the integrand factor $(t-\tau)$ with a positive exponent; here,
only a finite number if iterations are needed, since $\alpha>1$.
This  leads to
$\|u(t)\|_{\La}\le C_5+C_6\int_0^t \|u(\tau)\|_{\La}\,d\tau$ and finally to
$\|u(t)\|_{\La}\le C_5\exp(C_6t)$ by the classical Gronwall lemma.

To prove of \rf{d grad} under the stronger assumption $u_0\in\Ea$,
one could proceed in the same way, replacing the space $\La$
with $\Ea$ (and using the estimates of Lemma~\ref{LLL1}).
However, this argument would require additional
restrictions, such as inequalities of the form
$|f'(u)-f'(v)|\le C|u-v|(|u|^{q-2}+|v|^{q-2})$, which are not
fulfilled
for some  nonlinearities satisfying \rf{nla} with $q<2$.

Let us proceed in a slightly different way. First of all we have,
by \cite{BKW01,DroI06}, $\nabla u(t)\in
L^\infty([0,T],L^\infty(\R^d))$
for all $T>0$.
We rewrite the integral equation \rf{duham} in the following way
\begin{equation}
\label{int-duham}
\begin{split}
\nabla u(x,t)=&\nabla S_\alpha(t)u_0(x) \\
   &-\int_0^t\biggl(
	\int_{|y|\le |x|/2}+\int_{|y|\ge|x|/2}\biggr)
	\nabla p_\alpha(x-y,t-s) \nabla f(u(y,s))\,dy\,ds.
\end{split}
\end{equation}
It follows from condition~\rf{nla} that
$|f'(u)|\le C|u|^{q-1}$, hence, for every $u$ satisfying \rf{d sol}
we have
\begin{equation}
\label{nabf}
|\nabla  f(u(y,s))|\le C(1+|y|)^{-(q-1)(\alpha+d)}|\nabla u(y,s)|
\le C(1+|y|)^{-(q-1)(\alpha+d)},
\end{equation}
for a positive constant $C=C(T)$ and all $y\in\R^d$,
$s\in[0,T]$.
Combining \rf{nabf} with \rf{ST2} and with the decay estimate
$|\nabla p_\alpha(x,t)|\le Ct|x|^{-(\alpha+d+1)}$,
we get from~\rf{int-duham}
the preliminary inequality
\begin{equation}
\label{wee}
\begin{split}
|\nabla u(x,t)|\le C  & (1+|x|)^{-(\alpha+d+1)} +C(1+|x|)^{-
(\alpha+d+1)+q_1}
+C(1+|x|)^{-(q-1)(\alpha+d)}
\end{split}
\end{equation}
for some  constant $C=C(T)>0$, all $x\in\R^d$, $t\in
[0,T]$,
and with $q_1=d$.
Since $q>1$, now we can use this inequality to improve the estimate
in~\rf{nabf}.
This allows us to replace $q_1$ with some $0\le q_2<q_1$
and to improve also the estimate of the third term in~\rf{wee}.
After finitely many iterations of this argument
(more and more iterations are needed when $q$  approaches 1), we get
$|\nabla u(x,t)|\le C(T)(1+|x|)^{-(\alpha+d+1)}$ for all $x\in
\R^d$ and $t\in [0,T]$.
\endProof

\medskip

Let us now  recall a singular version of the Gronwall lemma.
This fact seems to be well-known,  we state it, however, in the form
which
is
the most suitable for our application and we prove it for the
completeness of the exposition.

\begin{lemma}\label{Gronwall}
Assume that a nonnegative and locally bounded function $h=h(t)$
satisfies the inequality
\begin{equation}
h(t)\leq C_1(1+t) +C_2 \int_0^t (t-\tau)^{-a} (1+\tau)^{-b}
h(\tau)\;d\tau\label{h est}
\end{equation}
for some $a\in(0,1)$, $b>0$, positive constants $C_1$ and $C_2$, and
all $t\geq 0$.
If $a+b>1$, then $h(t)\leq C(1+t)$ for all $t\geq 0$ and $C$
independent of $t$. The same conclusion holds true in the limit 
case $a+b=1$ under the weaker assumption
\begin{equation}
h(t)\leq C_1(1+t) +C_2 \int_0^t (t-\tau)^{-a} \tau^{-b}
h(\tau)\;d\tau\label{h est b}
\end{equation}
provided $C_2$ is sufficiently small.
\end{lemma}

\proof
If $a+b=1$, we deduce  from \rf{h est b}
the following inequality
$$
 h(t) \leq C_1(1+t) +C_2 K(a,b) \sup_{0\leq
\tau\leq t} h(\tau),
$$
where
$$
K(a,b)=
 \int_0^t (t-\tau)^{-a}
\tau^{-b}\;d\tau=\int_0^1 (1-s)^{-a} s^{-b}\;ds.  
$$
Consequently, $\sup_{0\leq
\tau\leq t} h(\tau) \leq  \frac{C_1}{1-C_2K(a,b)}(1+t)$ provided
$C_2<1/K(a,b)$.

In the case $a+b>1$, using \rf{h est}, we write $b=b_1+\eta$ with 
$a+b_1=1$ and
$\eta>0$,
and we fix $t_1>0$ such that
$C_2(1+t_1)^{-\eta}<1/K(a,b_1)$.
Now, splitting the integral in~\rf{h est} at $t_1$  yields
\begin{equation*}
h(t)\le C(1+t)+C_2K(a,b_1)(1+t_1)^{-\eta}
\sup_{0\le \tau\le t} h(\tau)
\end{equation*}
for some $C>0$ independent of $t$.
The conclusion of Lemma~\ref{Gronwall} now follows.
\endProof

\medskip
If the exponent $q$ in the assumptions on the nonlinearity \rf{hf1}
is larger than the time-critical value $\widetilde{q}$,
 we can improve the space decay estimates from
Proposition~\ref{prop1}
through the following space-time decay result.

\begin{theorem}\label{th up}
(i)
Let $\alpha\in (1,2)$. Assume that $u=u(x,t)$ is a solution of the
Cauchy problem \rf{eq}-\rf{ini}, where the nonlinearity $f$ satisfies
\rf{hf1} with $q>\widetilde{q}=1+(\alpha-1)/d$ and $u_0\in \La$. 
There exists $C>0$ (depending on $u_0$ but independent of $x,t$) such
that
\begin{equation}\label{up}
|u(x,t)|\leq Cp_\alpha(x,1+t) \quad \mbox{for all} \quad x\in \R^d
\quad\mbox{and}\;\;
t>0.
\end{equation}
The same conclusion holds true for $q=\widetilde{q}$ provided
$\|u_0\|_1$ is sufficiently small.

(ii)
Under the more stringent assumption $u_0\in\Ea$
we have also
\begin{equation}
\label{nup}
\|\nabla u(t)\|_{L^\infty_{\alpha+d+1}}
\le C (1+t).
\end{equation}
\end{theorem}

\proof
First recall that by estimates \rf{pbo} and \rf{Lp bound} with
$p=\infty$,
the solution satisfies
$$
\|u(t)\|_\infty \leq C(1+t)^{-d/\alpha}.
$$
Hence, to establish~\rf{up}, it suffices to prove
\begin{equation}
\|u(t)\|_\La\leq C(1+t).
\label{Lw est}
\end{equation}
Indeed,  the inequality
$$
g(x,t)\equiv \min \left\{(1+t)^{-d/\alpha}\,;\,
\frac{1+t}{(1+|x|)^{\alpha+d}}\right\}
\leq Cp_\alpha (x,t+1).
$$
is the consequence of the elementary estimate
$$
g(x,t) \leq (1+t)^{-d/\alpha} \min \left\{ 1\,;\,
 |x(1+t)^{-1/\alpha}|^{-\alpha-d}\right\}
$$
and   the asymptotic formula \rf{as P} (implying, in particular, that
$\min\{1\,;\,|x|^{-\alpha-d}\}\leq CP_\alpha(x)$ for all $x\in\R^d$
and a
constant $C>0$).

In the proof of \rf{Lw est}, we use the integral equation \rf{duham},
hence we
begin by the preliminary estimate (resulting from \rf{St 3} and from
the hypothesis \rf{hf1})
\begin{eqnarray*}
\|\nabla S_\alpha (t-\tau) f(u(\tau))\|_\La
&\leq& C(t-\tau)^{-1/\alpha}\|u(\tau)\|_\infty^{q-1}\|u(\tau)\|_\La\\
&&+
C(t-\tau)^{1-1/\alpha}\|u(\tau)\|_q^q.
\end{eqnarray*}
Moreover, since by \rf{pbo} and \rf{Lp bound} with $p=q$,
  the solution satisfies the decay estimate
\begin{equation}
\label{lqb}
\|u(\tau)\|^q_q\leq C(1+\tau)^{-d(q-1)/\alpha},
\end{equation}
we have the following inequalities
\begin{equation*}
\begin{split}
\int_0^t (t-\tau) ^{1-1/\alpha} \|u(\tau)\|_q^q\;d\tau
\leq
C\int_0^t (t-\tau) ^{1-1/\alpha}(1+\tau)^{-d(q-1)/\alpha} \;d\tau
\leq C(1+t)
\end{split}
\end{equation*}
which are valid for $1/\alpha+d(q-1)/\alpha\geq 1$.

Consequently,
after computing the $\La$-norm of equation \rf{duham}
and using  estimate \rf{St 2} we arrive at
\begin{equation}
\|u(t)\|_\La\leq C(1+t) + C\int_0^t (t-\tau) ^{-1/\alpha} (1+\tau)^{-
d(q-1)/\alpha} \|u(\tau)\|_\La\;d\tau.
\label{est sing}
\end{equation}

In the time-critical case $1/\alpha+d(q-1)/\alpha=1$ (i.e. for
$q=\widetilde{q}$)
 we
proceed
analogously, however, now we use the estimate
\begin{equation}
\label{lin}
\|u(\tau)\|_\infty \leq C \tau^{-d/\alpha}\|u_0\|_1
\end{equation}
with a constant $C$ independent of $u_0$ and $t$.
Hence, we obtain the following counterpart of inequality
\rf{est sing}
\begin{equation}
\|u(t)\|_\La\leq C(1+t) + C\|u_0\|_1^{q-1}
\int_0^t (t-\tau) ^{-1/\alpha} \tau^{
-
d(q-1)/\alpha} \|u(\tau)\|_\La\;d\tau.
\label{est sing bis}
\end{equation}

Finally, the singular Gronwall lemma (Lemma \ref{Gronwall})
applied
 to inequalities \rf{est sing} and \rf{est sing bis} completes
the proof of~\rf{up}.

To prove inequality \rf{nup} one should  follow exactly the same
argument as for the proof of~\rf{Lw est}, putting everywhere
$\Ea$-norms instead of the corresponding $\La$-norms,
and applying Lemma~\ref{LLL1} instead of Lemma~\ref{LL1}.
\endProof

\medskip

We conclude this section with estimates of self-similar solutions to
problem \rf{eq-sing}-\rf{ini-sing}.

\begin{corollary} \label{cor est U}
If the constant $M>0$ in \rf{ini-sing} is sufficiently small, 
then the corresponding solution of
problem \rf{eq-sing}-\rf{ini-sing} satisfies
\begin{equation}\label{self u p}
0\leq u_M(x,t)\leq Cp_\alpha(x,t), \qquad\hbox{for all $x\in\R^d$,
$t>0$},
\end{equation} 
with $C=C(M,\alpha,d)>0$ independent of $x$ and $t$.
\end{corollary}

\proof

Let us recall that the solution of \rf{eq-sing}-\rf{ini-sing} has been
constructed in \cite{BKW01}
as the limit of the rescaled functions
$u^\lambda(x,t)\equiv \lambda^d u(\lambda x, \lambda^\alpha t)$,
where $u=u(x,t)$ is the fixed solution of equation \rf{eq-sing}
supplemented
with the nonnegative initial datum $u(\cdot,0)=u_0\in
C^\infty_c(\R^d)$ such that
$\int_{\R^d} u_0(x)\;dx=M$. By Theorem \ref{th up}, used in the
critical case $q=\widetilde{q}$,
the rescaled family $u^\lambda$ satisfies
\begin{equation}\label{u l up}
|u^\lambda(x,t)|\leq C\lambda^d p_\alpha(\lambda x,
1+\lambda^{\alpha}t)
= Cp_\alpha(x,\lambda^{-\alpha} +t)
\end{equation}
for all $x\in \R^d$, $t>0$, and a constant  $C=C(M,\alpha,d)$
independent of
$x,t,\lambda$,
provided $M>0$ is sufficiently small.
Since $u^\lambda(x,t)\to U_M(x,t)$ as $\lambda\to\infty$  almost every
where in $(x,t)$
(see \cite[Lemma 3.7]{BKW01}), passing to the limit in \rf{u l up} we
complete the proof
of estimate~\rf{self u p}.
\endProof

%%%%%%%%%%%%%%%%%%%%%%%%%%%%%%%%%%%%%%%%%%%%%
\section{Asymptotic profiles}
\label{secP}

In this section, we derive the asymptotic expansions from Theorems
\ref{th21} and \ref{th22}.
Let us recall that
all positive constants, which appear here, are independent of $x$ and 
$t$ and are
denoted by the same letter~$C$.

\medskip

\noindent{\it Proof of Theorem \ref{th21}.}
Let us consider the nonlinear term appearing in the
integral equation~\rf{duham},
\begin{equation*}
{\mathcal{N}}(u)(t)\equiv\int_0^t\!\!\int \nabla_{\R^d}
 p_{\alpha}(x-y,t-s)
f(u(y,s))(s)\,ds.
\end{equation*}
In order to find  an asymptotics of $\mathcal N$ for large $|x|$,
we define two remainder functions $\RR(x,t)$ and
$\RR_1(x,t)$,
through the relations
\begin{equation}
\label{defR}
\begin{split}
\mathcal N(u)(x,t)
&=\int_0^t\!\!\int f(u(y,s)) \nabla p_\alpha(x,t-s)\,dy\,ds
+\RR_1(x,t)\\
&=-\frac{c_1\,x}{|x|^{\alpha+d+2}}
\int_0^t\!\! \int(t-s) f(u(y,s))\,dy\,ds
-\RR(x,t).
\end{split}
\end{equation}
Here, $c_1$ is the constant appearing in relation~\eqref{as nab P}.
Hence, it follows from the integral equation~\rf{duham}
that
\begin{equation}
\label{duhamR}
u(x,t)=S_\alpha(t)u_0(x)+\frac{c_1\,x}{|x|^{\alpha+d+2}}
\int_0^t\!\! \int_{\R^d} (t-s) f(u(y,s))\,dy\,ds
+\RR(x,t)
\end{equation}
and it remains
to estimate $\RR(x,t)$.

Computing the difference of the two expressions
of $\mathcal N$ in~\rf{defR} we deduce a bound  for  $\RR+\RR_1$,
implying
\begin{equation*}
\begin{split}
|\RR(x,t)| &\le |\RR_1(x,t)|\\
&+C|x|^{-(2\alpha+d+1)}
\int_0^t\!\!\int_{\R^d} \left[
\nabla p_\alpha(x,t-s)+\frac{c_1\,x}{|x|^{\alpha+d+2}}(t-s)
\right]|f(u(y,s))|\,dy\,ds.
\end{split}
\end{equation*}
Now, the asymptotic expansion~\rf{as nab P}, the assumption \rf{nla},
and the $L^q$-estimates \rf{lqb} lead to
\begin{equation}
\label{RR}
\begin{split}
|\RR(x,t)| &\le |\RR_1(x,t)|+C|x|^{-(2\alpha+d+1)}
\int_0^t\!\!\int_{\R^d} (t-s)^2 |f(u(y,s))|\,dy\,ds\\
&\le |\RR_1(x,t)|+C|x|^{-(2\alpha+d+1)} \,
t^2 \int_0^t(1+s)^{-(q-1)d/\alpha}\,ds.
\end{split}
\end{equation}

In order to estimate $\RR_1$,
we decompose it as $\RR_1=I_1+\cdots+I_4$,
where
\begin{equation*}
\label{inti}
\begin{split}
&I_1(x,t)\equiv
\int_0^t\!\!\int_{|y|\le |x|/2} [\nabla p_\alpha(x-y,t-s)-\nabla
p_\alpha(x,t-s)]
\cdot f(u(y,s))\,dy\,ds,\\
&I_2(x,t)\equiv
-\int_0^t \biggl(\int_{|y|\ge |x|/2}f(u(y,s))\,dy\biggr)
	\nabla p_\alpha(x,t-s)\,ds,\\
&I_3(x,t)\equiv
\int_0^t\!\!\int_{|y|\ge |x|/2,\; |x-y|\ge |x|/2}
\nabla p_\alpha(x-y,t-s) f(u(y,s))\,dy\,ds,\\
&I_4(x,t)\equiv
\int_0^t\!\!\int_{|y|\le |x|/2} \nabla p_\alpha(y,t-s)
f(u(x-y,s))\,dy\,ds.
\end{split}
\end{equation*}

In our next two estimates, we use the
inequality (which is a consequence of the
$L^\infty$-bound of the solution, see~\eqref{Lp bound})
\begin{equation}
\label{peu}
|u(y,s)|^q\le
C(1+|y|)^{-(\alpha+d)}(1+s)^{-(q-1)d/\alpha}\|u(s)\|_\La .
\end{equation}
This leads to
\begin{equation}
\label{iiib}
\begin{split}
 |I_1(x,t)| &
   \le C|x|^{-(\alpha+d+2)}
    \int_0^t (t-s) \!\!\int_{|y|\le |x|/2}
     |y|\, |u(y,s)|^q\,dy\,ds\\
    &\le C|x|^{-(\alpha+d+2)}\,t
      \int_0^t (1+s)^{-(q-1)d/\alpha} \| u(s)\|_\La \,ds.
      \end{split}
      \end{equation}
Here, we have applied also
the Taylor formula and the bound~\rf{n2p}.

The next two integrals can be bounded by the same quantity, indeed
\begin{equation}
\label{iiib2}
\begin{split}
 |I_2(x,t)|+|I_3(x,t)| &\le C|x|^{-(\alpha+d+1)}
     \int_0^t (t-s)\int_{|y|\ge |x|/2} |u(y,s)|^q\,dy\,ds\\
     & \le C|x|^{-(2\alpha+d+1)}\,t \int_0^t (1+s)^{-(q-1)d/\alpha}
\|u(s)\|_\La\,ds.
\end{split}
\end{equation}

The  estimate for the last term is
 \begin{equation}
  \label{i4}
  |I_4(x,t)|
     \le C|x|^{-q(\alpha+d)}\int_0^t
(t-s)^{-1/\alpha}\|u(s)\|_\La^q\,ds.
\end{equation}
Since we are assuming $\alpha>1$, when we compare the exponents
of $|x|$ in inequalities~\rf{RR} and~\rf{iiib}-\rf{iiib2}, we see that
\begin{equation}
\label{rrf}
|\RR(x,t)|\le C\, \max\bigl\{ |x|^{-q(\alpha+d)} \,;\,
|x|^{-(\alpha+d+2)} \bigr\}
\quad\mbox{for all $|x|\ge1$ and  $t\in (0,T]$},
\end{equation}
where $C=C(T)>0$ is uniformly bounded with respect to $x$ and $t$,
in any time interval $t\in[0,T]$.
Part (i) of Theorem~\ref{th21} now follows.

To establish Part (ii),
we have only to improve the estimate of the integral~\rf{i4}.
We can do it using, in a slightly deeper way,
the properties of the fundamental solution $p_\alpha(x,t)$.
In particular, its radial symmetry implies that
$$\int_{|y|\le |x|/2} \nabla p_\alpha(y,t-s)\,ds=0,$$
so that
\begin{equation}
\label{newi4}
I_{4}(x,t)\equiv
\int_0^t\!\!\int_{|y|\le |x|/2} \nabla p_\alpha(y,t-s)\cdot
	\bigl[f(u(x-y,s))-f(u(x,s)) \bigr]\,dy\,ds.
\end{equation}
Owing to the more stringent assumption
$u_0\in \Ea$ and by Proposition~\ref{prop1},
we deduce from the mean value theorem applied to $f(u)$
(recall that $|f'(u)|\le C|u|^{q-1}$)
\begin{equation}
\label{ini41}
|I_{4}|\le C|x|^{-q(\alpha+d)-1}\int_0^t\|u(s)\|_\La^{q-1}\|\nabla
u(s)\|_{L^\infty_{\alpha+d+1}}
\,ds.
\end{equation}
Replacing inequality~\rf{i4} with this new estimate shows that
the bound~\rf{rrf} of the remainder term
can be improved into
\begin{equation}
\label{rrfim}
|\RR(x,t)|\le C\, \max\bigl\{ |x|^{-q(\alpha+d)-1} \,;\,
|x|^{-(\alpha+d+2)} \bigr\}
\quad\mbox{for all $|x|\ge1$ and $t\in (0,T]$}.
\end{equation}
Hence, Part (ii) of Theorem \ref{th21} follows.

Let us prove assertion (iii).
When the solution satisfies
the additional estimate~\rf{u p a}
(recall that, by Theorem~\ref{th up}, such an estimate holds true at 
least when either
$q>\widetilde{q}$ or  $q=\widetilde{q}$ and $\|u_0\|_1$
is small enough),
 we have  $\|u(t)\|_\La\le C(1+t).$
In this case, it is easy to construct an exponent $N=N(\alpha,d,q)$
such that
\begin{equation}
\label{rrf2}
|\RR(x,t)|\le C(1+t)^N
\max\bigl\{ |x|^{-q(\alpha+d)} \,;\, |x|^{-(\alpha+d+2)} \bigr\}
\quad\hbox{for all $|x|\ge1$, $t>0$}.
\end{equation}
Let us explain why   $N\le 3$.
It follows directly from~\rf{RR} and from \rf{iiib}-\rf{i4}
that $N\le \max\{3\,;\,q+1-1/\alpha\}$.
However, if $q> 2+1/\alpha$, then
 we can replace
estimate~\rf{i4} with
\begin{equation}
\label{i42}
\begin{split}
|I_4(x,t)| &\le C|x|^{-(\alpha+d+2)}\int_0^t (t-s)^{-1/\alpha}
\|u(s)\|_\infty^{q-1-2/(\alpha+d)}\|u(s)\|_\La^{1+2/(\alpha+d)} \,ds\\
   &\le C(1+t)^3 |x|^{-(\alpha+d+2)}.
\end{split}
\end{equation}

If, moreover,  $\nabla u$
satisfies the additional pointwise estimate~\rf{nupi}
then we can precise in a similar way the bound~\rf{rrfim}.
Namely, we can replace $C=C(T)$ in \rf{rrfim} with $C(1+t)^3$.
Next, the proof of this claim relies either on inequality 
\eqref{ini41} 
if
$1<q\le2$
or on the following new estimate of $I_4$ when $q>2$
$$ |I_{4}(x,t)|\le C|x|^{-(\alpha+d+2)}\int_0^t 
\|u(s)\|_\infty^{q-1-1/(\alpha+d)}\|u(s)\|_\La^{1/(\alpha+d)}\|\nabla
u(s)\|_{L^\infty_{\alpha+d+1}}
  \,ds.$$
The estimates of the other terms remain unchanged.
The proof of Theorem~\rf{th21} is now complete.
\endProof

%%%%%
\medskip

\noindent{\it Proof of Theorem~\ref{th22}.}
Let $u_M$ be a self-similar solution of~\rf{eq-sing}-\rf{ini-sing},
satisfying estimate~\rf{self u p 0}.
We consider the integrals $I_1$, $I_2$, $I_3$ and $I_{4}$
and also the remainder term $\RR$, defined as in the proof of Part~(i)
of Theorem~\ref{th21}.
We treat all these terms proceeding as before,
but replacing everywhere estimate~\rf{peu}
with the estimate (deduced from~\rf{self u p 0})
\begin{equation}
\label{peus}
|u_M(y,s)|^q \le Cs^{-dq/\alpha} P_\alpha^q(y/s^{1/\alpha}),
\end{equation}
with $q=\widetilde{q}$ and $C=C(M)$,
then making the change of variables $y\mapsto ys^{1/\alpha}$ in all
the
space integrals.
After some simple computations, we arrive at
\begin{equation*}
|\RR(x,t)| \le Ct^{-d/\alpha} \biggl[
\Bigl( \frac{|x|}{t^{1/\alpha}} \Bigr)^{-(2\alpha+d+1)} +
\Bigl( \frac{|x|}{t^{1/\alpha}} \Bigr)^{-(\alpha+d+2)} +
\Bigl( \frac{|x|}{t^{1/\alpha}} \Bigr)^{-\widetilde{q}(\alpha+d)}
\biggr].
\end{equation*}
Recalling that $f(u)=b u^{\widetilde{q}}$,
applying~\rf{duhamR} to $u_M$ we get
\begin{equation*}
u_M(x,t)= Mp_\alpha(x,t)+t^{1+1/\alpha}\cdot
 \frac{c_1\alpha^2}{\alpha+1}
\biggl(\int U_M(y)^{\widetilde{q}}\, dy\biggr)
\frac{b\cdot x}{|x|^{\alpha+d+2}}
+\RR(x,t).
\end{equation*}
Now, passing to self-similar variables, we deduce that, for all
$x\in\R^d$,
\begin{equation*}
U_M(x)=M P_\alpha(x) +\frac{c_1\alpha^2}{\alpha+1}
\|U_M\|_{\widetilde{q}}^{\widetilde{q}} \,
\frac{b\cdot x}{|x|^{\alpha+d+2}}
+\RR_M(x),
\end{equation*}
where
\begin{equation*}
\RR_M(x)=O\Bigl( \max\bigl\{ |x|^{-\widetilde{q}(\alpha+d)} \,;\,
	|x|^{-(\alpha+d+2)} \bigr\}
	\Bigr),\qquad\hbox{as $|x|\to\infty$}.
\end{equation*}
Theorem~\rf{th22} is now established.
\endProof

\begin{remark}
We conclude observing that
the above  expression of the remainder term $\RR_M(x)$ can be 
simplified
distinguishing the
 two cases $d=1$ and $d\ge2$.
Indeed,  an elementary
calculation shows that
\begin{itemize}
\item[{1.}]
In the one dimensional case $d=1$ (hence, $\widetilde{q}=\alpha$,
and the assumption $\widetilde{q}>q^*$ reads $\alpha>\sqrt 2$), we
have
\begin{equation*}
\RR_M(x)=
\left\{
\begin{array}{lcc}
O\bigl( |x|^{-\alpha(\alpha+1)} \bigr) &\mbox{if}&
\sqrt 2<\alpha\le \sqrt 3,\\
O\bigl( |x|^{-(\alpha+3)} \bigr)& \mbox{if}&
 \sqrt
3\le \alpha<2,
\end{array}
\right.
\quad \mbox{as}\quad |x|\to \infty.
\end{equation*}

\item[{2.}]
For $d\ge2$, it follows
\begin{equation*}
\RR_M(x)=O\bigl( |x|^{- \widetilde{q}(\alpha+d)} \bigr)
\quad\mbox{as}\quad |x|\to\infty.
\end{equation*}
\end{itemize}
\end{remark}

\begin{remark}
Analogously, as in Theorem \ref{th21}, one could remove the
restriction $\widetilde{q}>q^*$ from Theorem \ref{th22}, provided we 
have the additional weighted estimate 
\begin{equation}
\label{nupis}
\|\nabla u_M(t)\|_{L^\infty_{\alpha+d+1}}
\le C t.
\end{equation}
We expect that inequality \rf{nupis} can be proved using the 
scaling argument from the proof of Corollary \ref{cor est U}, below.
This reasoning would require, however, some improvements of 
estimates from \cite{BKW01}. 
We skip other details because the goal of this work was to present a 
method of deriving asymptotic expansions of solutions rather than to 
study the most general case.
\end{remark}

%%%%%%%%%%%%%%%%%%%

\bibliographystyle{amsplain}

\end{document}